\documentclass[12pt,reqno]{amsart}%[11pt]{article}

\usepackage{showkeys}
\usepackage{amsmath,enumitem,empheq}
\usepackage{latexsym}
\usepackage{amssymb}
\usepackage{amsfonts}
\usepackage[usenames]{color}

\usepackage{graphics}
\include{PDF}
\renewcommand{\proof}{\par\noindent{\it Proof.\ \ }}
\def\qed{\ifmmode\square\else\nolinebreak\hfill
$\square$\fi\par\vskip12pt}

\def\ov{\overline} 
\def\l{\langle} \def\r{\rangle}

 \def\ZZ{\mathbb Z}
\def\BB{{\mathcal B}}

\def\Cay{{\sf Cay}}

\def\D{{\rm D}} 
\def\S{{\rm S}} 
 \def\M{{\rm M}}

\def\C{{\bf C}} \def\N{{\bf N}} 
\def\Z{{\bf Z}} 

\def\mod{{\sf mod~}} 
\def\char{{\sf \,char\,}}

\def\Aut{{\sf Aut}}  
\def\K{{\sf K}}

 \def\AS{{\rm AS}}

\def\Ome{{\it\Omega}}
\def\Ga{{\it\Gamma}} \def\Sig{{\it \Sigma}} \def\Del{{\it \Delta}}

\def\a{\alpha} \def\b{\beta} \def\g{\gamma} \def\s{\sigma}
\def\t{\tau}

\def\GammaL{{\rm \Gamma L}}

\def\A{{\rm A}}

\def\PSL{{\rm PSL}}\def\PGL{{\rm PGL}}
\def\GL{{\rm GL}} 

 \def\PSU{{\rm PSU}}

  \def\D{{\rm D}}

\newtheorem{theorem}{Theorem}[section]%
\newtheorem{lemma}[theorem]{Lemma}%
\newtheorem{proposition}[theorem]{Proposition}%

\begin{document}

\title[Metacirculants]
{A classification of tetravalent edge-transitive   metacirculants of odd order}
%\thanks{This work forms a part of an ARC grant project.}

%\address{School of Mathematics and Statistics\\
%The University of Western Australia\\
%Crawley, WA 6009, Australia}
%

\author[Song]{Shu Jiao Song}
\address{School of Mathematics and Statistics\\
The University of Western Australia\\
Crawley, WA 6009, Australia}

\email{shu-jiao.song@uwa.edu.au}
%\email{shu-jiao.song@uwa.edu.au}

\date\today%{1/9/2002}%{21/8/2002}%{10/5/2002}

\begin{abstract}
 In this paper a classification of tetravalent edge-transitive metacirculants is given. It is shown that a tetravalent  edge-transitive  metacirculant $\Ga$ is a normal graph except for four known graphs. If further, $\Ga$ is a Cayley graph of a non-abelian metacyclic group, then $\Ga$ is half-transitive.
\end{abstract}

\maketitle

% \footnotetext{* This work was partially supported by an NNSF(K1020261) and an ARC Discovery Grant.}
%\footnotetext{$\star$ Corresponding author: Department of
%Mathematics, Yunnan University, Kunming Yunnan, 650091, P. R.
%China,}
%\author{ \qquad Xiao-hui Zhang  \\ {\small{{ Department
%of Mathematics, Yunnan University, Kunming, Yunnan}}},
% \\{\small{ 650091, P. R. China }}}

%A group $G$ is called {\it metacyclic} if there is a normal cyclic subgroup $N$ such that 
%the factor group $G/N$ is cyclic.
%Thus $G=\ZZ_n.\ZZ_m$, and $G=\l a,b\r$ for some elements $a,b$ such that $\l a\r$ is normal in $G$.
%If the extension $\ZZ_n.\ZZ_m$ is split, then $G$ is called a {\it split metacyclic group}.
%In this case, $G=\l a\r{:}\l b\r$.

\section{Introduction}
Throughout this paper graphs are assumed to be finite,  simple, connected and undirected, unless stated otherwise. 

A group is called {\it metacyclic} if and only if there exists a cyclic normal subgroup $K$ of $G$ such that  $G/K$ is cyclic. A metacyclic group is called  {\it split} if it is a split extension of a cyclic group by a cyclic group, namely, $G=H{:}K$ is a semidirect product of a cyclic subgroup $H$ by a cyclic subgroup $K$. A graph is called a {\it metacirculant} if it contains a metacyclic vertex-transitive automorphism group,
introduced by Maru\v si\v c and \'Sparl \cite{M-Sparl}. In this paper, a Cayley graph of  a split metacyclic group $G$ is  called  a {\it split metacirculant} of $G$. 

A graph $\Ga=(V,E)$ is called $X$-{\it edge-transitive}
if a subgroup $X \leq \Aut\Ga$ is transitive on $E$. 
%It is called $X$-{\it vertex-transitive} if $X$ is transitive on $V$. An edge-transitive metacirculant graph is both vertex-transitive and
%edge-transitive. Certain special classes edge-transitive metacirculants have ben investigated, see
%\cite{ACMX,K} for a classification of edge-transitive circulants.
%
An {\it arc} of a graph is an ordered pair of
adjacent vertices. A graph $\Ga=(V,E)$ is called $X$-{\it
half-transitive} if $X\le\Aut\Ga$ is transitive on both $V$ and $E$,
but not transitive on the arc set of $\Ga$; in particular, if
$X=\Aut\Ga$ then $\Ga$ is simply called a {\it half-transitive
graph}.

%It is a generalisation of the concept of {\it matecirculant} (introduced in \cite{AP}), 
%which has a metacyclic vertex-transitive automorphism group with certain restrictive condition.
%It is shown in \cite[Thm\,1.3]{LSW} that the class of   metacirculants is indeed larger than that of metacirculant. In this paper, we simply called   metacirculant graphs 

In the literature, the class of metacirculants provides a rich source of various research projects, see for instance \cite{AMN,MM,San,S,T}. It has been therefore studied extensively in the past 30 years. Special classes of metacirculants have been well investigated, see \cite{ACMX,DMM,K,Li,LPSW,LSW-1, M} for edge-transitive Cayley graphs of special metacyclic groups; \cite{Li-Sims,Xu-half} for half-transitive metacirculants of special order;   \cite{LSW} for cubic metacirculants; \cite{M-Sparl,Sparl-meta,Sparl-meta2} for tetravalent half-transitive metacirculants.
%The investigation of half-transitive graphs was initiated with a
%question of Tutte \cite{Tutte} regarding the existence for even
%valencies $2k\geq 4$.
%In 1970, Bouwer \cite{I.z} constructed
%the first family of half-transitive graphs, and since then,
%constructing and characterizing half-transitive graphs  has been an active
%topic in algebraic graph theory (see, for example,
%\cite{A&X,Marusic,M-N}). Many interesting families of
%half-transitive graphs are constructed as metacirculants, see
%\cite{Li-Sims,Xu-half} for a family of half-transitive
%metacirculants of prime-power order. 
However, despite all the efforts there are still numerous questions about these graphs that need to be answered. Indeed, even the family of tetravalent half-transitive metacirculants, which has the smallest admissible valency for a half-transitive graph, seems to be too difficult to be completely classified.
%was studied the most, seems to be too rich to be completely classified, see \cite{FKXZ,FKXZ-2,FWZ,LLL,LLZ,Mar-4,PP,Sparl-4,ZF} for the study of tetravalent half-transitive graphs, and see \cite { M-Sparl,Sparl-meta,Sparl-meta2} for the study of tetravalent half-transitive metaciruclants.
 In this paper, we  study tetravalent edge-transitive metacirculants of odd order and give a complete classification of these graphs.

A graph $\Ga$ is a {\it Cayley graph} if there exists a group $G$ and a
subset $S\subset G$ with $S=S^{-1}:=\{s^{-1} \mid s\in S\}$ and $1\not\in S$ such
that the vertices of $\Ga$ may be identified with the elements of
$G$ in the way that $x$ is adjacent to $y$ if and only if
$yx^{-1}\in  S$, such a graph is denoted by
$\Cay(G,S)$. A Cayley graph $\Ga = \Cay(G, S)$ has an automorphism
group
\[\hat G = \{\hat g: x\mapsto xg, \mbox{for all } x\in G \mid g\in G\},\]
consisting of right multiplications by elements $g \in G$.
If further $\hat G$ is normal in $\Aut\Ga$,
then $\Ga$ is called a {\it normal Cayley} graph.

The main result of this  paper is stated in the following theorem.

\begin{theorem}\label{thm-1}
Let $\Ga$ be a connected edge-transitive tetravalent   metacirculant. 
Then  $\Ga$ is   split metacirculant of $G$ and one of the following holds:
 \begin{itemize}
 \item[1)] $\Ga$ is an arc-transitive normal Cayley graph of an abelian group such that $\Ga\neq \K_5$;
% \item[2)] $\Ga=\Cay(G,S)$ is a normal Cayley graph of a abelian group $G$ such that  $S=\{x,y,x^{-1},y^{-1}\}$  with $o(x)=o(y)$ and $\l x,y\r=G$, and $\Aut\Ga=G{:}(\Aut\Ga)_\a$ where $\ZZ_2<(\Aut\Ga)_\a\leqslant \D_8$;

 \item[2)] $\Ga$ is a half-transitive normal Cayley graph of a non-abelian split metacyclic group $G$ where $G\neq  \ZZ_{7}{:}\ZZ_3$, $\ZZ_{11}{:}\ZZ_5$ or $\ZZ_{23}{:}\ZZ_{11}$. The automorphism group $\Aut\Ga=G{:}\ZZ_2$.
  
 \item[3)] $\Ga=\Cay(G,S)$ with $G=\ZZ_5$, $\ZZ_{7}{:}\ZZ_3$, $\ZZ_{11}{:}\ZZ_5$ or $\ZZ_{23}{:}\ZZ_{11}$ and
 $\Ga$ is one of $\ n$ non-isomorphic edge-transitive graphs. Among them, one is $s$-arc-transitive with $(\Aut\Ga,(\Aut\Ga_\a), s)$ listed in Table 1, and the other $n-1$ are half-transitive with $\Aut\Ga=G{:}\ZZ_2$.
% \begin{itemize}
% \item[i)] $G=\ZZ_5$, $\Aut\Ga=\S_5$, $\Ga$ is $2$-arc-transitive graphs;
% \item[ii)] $G=\ZZ_{7}{:}\ZZ_3$, $\Ga$ is one of  three non-isomorphic edge-transitive graphs. Among them, one is arc-transitive with $\Aut\Ga=\PGL(2,7)$ and vertex stabiliser is $D_{16}$, the other two are half-transitive graphs with vertex stabiliser is $\ZZ_2$;
% \item[iii)] $G=\ZZ_{11}{:}\ZZ_5$, $\Ga$ is one of  six non-isomorphic edge-transitive graphs. Among them, one is 2-arc-transitive with $\Aut\Ga=\PGL(2,11)$ and vertex stabiliser is $\S_4$, the other five are half-transitive graphs with vertex stabiliser is $\ZZ_2$;
%\item[iv)] $G=\ZZ_{23}{:}\ZZ_{11}$, $\Ga$ is one of  eleven non-isomorphic edge-transitive graphs. Among them, one is arc-transitive with $\Aut\Ga=\PGL(2,23)$ and vertex stabiliser is $\S_4$, the other ten are half-transitive graphs with vertex stabiliser is $\ZZ_2$;
%
% \end{itemize}

\[\begin{array}{l l l l l l l l l} \hline
\Aut\Ga& &(\Aut\Ga)_\a & & s & & n & &G \\ \hline
\S_5  &\qquad & \S_4 & \qquad& 2 & \qquad& 1 &\qquad & \ZZ_5\\
\PGL(2,7)                      & & \D_{16}      & & 1 & \qquad& 3 &\qquad & \ZZ_7{:}\ZZ_3\\
\PGL(2,11)   & &  \S_4 & & 2 &\qquad& 6 &\qquad  & \ZZ_{11}{:}\ZZ_5\\
\PSL(2,23)                     & & \S_4          & & 2 & \qquad& 11 &\qquad & \ZZ_{23}{:}\ZZ_{11}\\ \hline

\end{array}\]

\centerline{{\bf Table 1}\ }

\end{itemize}

\end{theorem}

A metacirculant $\Ga$ is called a {\it Sylow-circulant} if $\Aut\Ga$ has a transitive metacyclic subgroup of which all Sylow subgroups are cyclic. 

{\bf Hypothesis (*):} Let $G$ be a split metacyclic group such that all Sylow subgroups of $G$ are cyclic.
Then $G$ can be written as \[G=H\times K=\l c\r\times \l a,b\r=\l a\r{:}\l bc\r,\] where \[K=\l a,b\mid a^m=b^n=1,a^b=a^r\r, r^n\equiv 1\pmod m,\] such that no Sylow $p$-subgroup of $\l b\r$ lies in the centre of $G$.
Let $n_0$ be the smallest integer such that $\l b^{n_0}\r\in \Z(G)$.

We remark that $H$ may be a trivial subgroup of $G$. In this case, $c$ is the identity of $G$. Recall that the Euler Phi-function $\phi(m)$ is the number of positive integers that are less than m and coprime to $m$.
\begin{theorem}\label{thm-2}
Let $G=H\times K=\l c\r\times\l a,b\r$ and $n_0$ be as hypothesis (*). Let $\Ga=\Cay(G,S)$ be a connected  tetravalent edge-transitive graph. Then \[S=\{cb^j,c^{-1}ab^{j},c^{-1}b^{-j},c(ab^j)^{-1}\},\]where $1\leqslant j<n_0$ and $(j,n)=1.$

In particular, there are exactly $\phi(n_0)/2$ non-isomorphic tetravalent  edge-transitive Cayley graphs of $G$, except for $G=\ZZ_{11}{:}\ZZ_5$ and $\ZZ_{23}{:}\ZZ_{11}$ which has $3$ and $6$ non-isomorphic tetravalent  edge-transitive Cayley graphs of $G$, respectively.
\end{theorem}

%\begin{corollary}
%Let $G=\l a\r{:}\l b\r=\ZZ_m{:}\ZZ_n$ such that $\Z(G)=1$.
% Let $\Ga=\Cay(G,S)$ be a connected edge-transitive tetravalent graph. 
% Then $S=\{b^j,ab^j,b^{-j},(ab^j)^{-1}\}$ for some integer $j$ coprime to n.
% 
% In particular, there are exactly $\phi(n)/2$ non-isomorphic edge-transitive tetravalent Cayley graphs of $G$  for $G\neq \ZZ_{11}{:}\ZZ_5$.
% 
% If $G=\ZZ_{11}{:}\ZZ_5$, then there are 3 non-isomorphic edge-transitive tetravalent Cayley graphs of $G$.
% 
%\end{corollary}

\section{Notation and preliminaries}

Let $\Ga=\Cay(G,S)$ be a Cayley graph of a group $G$. Then  $\Aut\Ga$
contains the a regular subgroup which is isomorphic to $G$. For an automorphism group $X$ such that $G\le X\le\Aut\Ga$, we have
$X= GX_\a$, where $\a\in V$.
Conversely, by \cite[Proposition 4.3]{Xu}, if $\Aut\Ga$ has a regular subgroup $G$, then $\Ga\cong\Cay(G,S)$ for some $S\subset G.$

Besides the regular automorphism subgroup $G$, $\Aut\Ga$ has another important subgroup:
\[\Aut(G,S)=\{\s\in\Aut(G)\mid S^\s=S\}.\]
Let $\a$ be the vertex of $\Ga$ corresponding to the identity of
$G$. Then $S=\Ga(\a)$,  $\Aut(G,S)<\Aut\Ga$ and $\Aut(G,S)$ fixes $\a$ and
normalises $ G$. Moreover, by \cite[Lemma 2.1]{Godsil}, we have the
following lemma.

\begin{lemma}\label{N(G)}
The normalizer $\N_{\Aut\Ga}( G)=G{:}\Aut(G,S)$.
\end{lemma}

For any automorphism $\s\in\Aut(G)$, the Cayley graphs $\Cay(G,S)$
and $\Cay(G,S^\s)$ are isomorphic. A Cayley graph $\Cay(G,S)$ is
called a {\it CI-graph} (CI stands for {\it Cayley Isomorphism}) if,
whenever $\Cay(G,T)\cong\Cay(G,S)$, we have $T=S^\s$ for some
automorphism $\s\in\Aut(G)$.

\begin{proposition}\label{CI-criterion} {\rm(\cite[Theorem~4.1]{val<q})}
Let $\Ga$ be a Cayley graph of a group $G$, and let $A=\Aut\Ga$.
If $G$ is of odd order and $(|G|,|A_\a|)=1$, then $\Ga$ is
a CI-graph of $G$.
\end{proposition}

Let $G$ be a transitive permutation group on $\Ome$. Then the $G$-orbits in $\Ome\times\Ome$ are called the $G$-{\it orbitals} on $\Ome$. 
For each orbital $\Del$, there is a {\it paired orbital }$\Del^*$, namely $\Del^*:=\{(\b,\a)\mid (\a,\b)\in \Del\}$.
If $\Del=\Del^*$, we say that the orbital $\Del$ is self-paired. 

For a  $G$-orbital  $\Del:=(\a,\b)^G$, there is a di-graph $\Ga$ with vertex set $\Ome$, and arc set $(\a,\b)^G$. $\Ga$ is called an {\it orbital graph} of $G$ on $\Ome$. An orbital graph is undirected if and only if $\Del$ is self-paired.

\vskip0.1in

 Let $\Ga=(V,E)$ be a connected $X$-edge-transitive graph. Let $N$
be an intransitive normal subgroup of $X$. Let $\BB$ be the set of $N$-orbits on
$V$, which is sometimes denoted by $V_N$. The {\it quotient graph}
$\Ga_\BB$ is the graph with vertex set $\BB$ such that two
vertices $B,B'\in\BB$ are adjacent in $\Ga_\BB$ if and only
if there exist $\a\in B$ and $\a'\in B'$ with $\{\a,\a'\}\in E$.
The graph $\Ga_\BB$ is sometimes denoted by $\Ga_N$, called the {\it
normal quotient} of $\Ga$ induced by $N$. Let $K$ be the kernel of $X$ acting
on $\BB$. Then $X/K\lesssim \Aut\Ga_N$ and $X/K$ is edge-transitive on $\Ga_N$.

As usual, for a set $\pi$ of primes, denote by $\pi'$ the set of
primes which are not in $\pi$.
For a group $L$, a {\it Hall $\pi$-subgroup} $H$ of $L$  is a subgroup
such that each prime divisor of $|H|$ lies in $\pi$ and the
index $|L:H|$ is coprime to all primes in $\pi$.
Let $n$ be an integer, we denote $n_{\pi}$ the divisor of $n$ such that  each prime divisor of $n_{\pi}$ lies in $\pi$ and $n/n_{\pi}$ is coprime to $n_{\pi}.$

\begin{lemma}\label{Hall}
Let $X$ be an insoluble group which has Hall $\pi$-subgroups. Then each normal or subnormal subgroup of $X$ has Hall $\pi$-subgroups.
\end{lemma}
\proof
Let $T\lhd M\lhd X$. Since $M\lhd X$, we have $MX_{\pi}$ is a subgroup of $X$. So \[|X|\geqslant |MX_{\pi}|=\frac{|X_{\pi}||M|}{|M\cap X_{\pi}|}.\]
Hence $|X|_{\pi}\geqslant |X|_{\pi}\cdot\frac{|M|_{\pi}}{|M\cap X_{\pi}|}. $ It follows that $|M|_{\pi}=|M\cap X_{\pi}|$. So $M$ has  Hall $\pi$-subgroups. Similarly, $T$ has Hall-$\pi$ subgroups.\qed 

Let $G$ be a  metacyclic group. Then any subgroup or quotient group of $G$ is a metacyclic group. The following lemma shows that $G$ can not be written as a direct product of two isomorphic non-cyclic groups.

\begin{lemma}\label{non-iso-direct}
Let $G=G_1\times G_2$ be a metacyclic group such that $G_1\cong G_2$. Then $G_1\cong G_2$ are cyclic groups
\end{lemma}

\proof
Suppose $G_1=H_1. K_1\cong G_2=H_2. K_2$, and $H.K=G=G_1\times G_2$ be a metacyclic group, where $H_1\cong H_2$, $K_1\cong K_2$ and $H,K$ are cyclic groups. 
We first assume that $G_1\cong G_2$ are non-abelian groups.
Since $H_i\times \C_{K_i}(H_i)\lhd G$ $(i=1,2)$, we have $G_1/(H_1\times \C_{K_1}(H_1))\times G_2/(H_2\times \C_{K_2}(H_2))\cong (G_1\times G_2)/\Big(H_1\times \C_{K_1}(H_1)\Big)\times \Big(H_1\times \C_{K_1}(H_1)\Big)$
Since $G_i$ are not abelian,  there is a prime $p$ such that the Sylow $p$-subgroup $(K_i)_p$ of $K_i$ is not  normal in $G_i$ with $i=1,2$. So $p$ divides $G_i/\Big(H_i\times \C_{K_i}(H_i)\Big)$. Thus $\ZZ_p^2\lesssim G_1/(H_1\times \C_{K_1}(H_1))\times G_2/(H_2\times \C_{K_2}(H_2))$ is not cyclic. On the other hand, each normal subgroup of $G$ is contained in $H_1\times \C_{K_1}(H_1)\times H_2\times \C_{K_2}(H_2)$, so does $H$. Thus $G/\Big(H_1\times \C_{K_1}(H_1)\Big)\times \Big(H_1\times \C_{K_1}(H_1)\Big)\lesssim G/H$ which is not possible as $G/H$ is cyclic. So $G_i$ $(i=1,2)$ are abelian groups.

%Suppose that $P_i:=(G_i)_p$ is a non-cyclic Sylow p-subgroup of $G_i$ $(i=1,2)$. 
%We have $P:=P_1\times P_2\leqslant G$. Let $\phi(P_i)$ be the Frattini subgroup of  $P_i$, that is, the intersection of all maximal subgroups of $P_i$. Then $P_i/\phi(P_i)$ is elementary abelian $p$-group. Since $P_i$ is non-abelian metacyclic group, we have $P_1/\phi(P_1)\cong \ZZ_p^2$. On the other hand, $\phi(P_i)\char P_i\lhd P$. Hence $\ZZ_p^4\cong P_1/\phi(P_1)\times P_2/\phi(P_2)\cong P/(\phi(P_1)\times \phi(P_2))$ which is not possible as $P$ is metacyclic group and so is $P/(\phi(P_1)\times \phi(P_2)) $ .
%
%Now suppose each Sylow subgroup of $G_i$ is cyclic. Then each normal subgroup of $G$ is contained in $H_1\times \C_{K_1}(H_1)\times H_2\times \C_{K_2}(H_2).$ While $G_1$
%
%
%Then there is a prime $p$ such that the Sylow $p$-subgroup $(K_i)_p$ of $K_i$ is not  normal in $G_i$ with $i=1,2$. Thus $(K_i)_p$ is not contained in $H$.
%Thus $p$ divides $|K_i/H\cap K_i|\leq |G_i/H\cap G_i|$.  Since $G=H{:}K=G_1\times G_2$, we have $K=(G_1\times G_2)/H$Thus $\ZZ_p^2\leqslant K_1/(H\cap K_1)\times K_2/ (H\cap K_1) \lesssim G/H$ which is not possible. So $G_i$  ($i=1,2$) are abelian groups.
%

Now suppose $G_1\cong G_2$ are non-cyclic abelian groups.
Then there is a prime $p$ such that $p\mid (|H_1|,|K_1|)$.
So $\ZZ_p^4\leqslant G$ which is not possible as $G$ is a metacyclic group. So $G_i$ $(i=1,2)$ are cyclic groups.\qed

%Let $H$ be a maximal normal cyclic subgroup of $G$.
%Then $G/H$ is a cyclic group.
%Moreover, \[H\cap K_1=K_1\cap \Z(G)\cong K_2\cap \Z(G).\]
%As $(K_i)_p\not\leqslant \Z(G_i)$ with $i=1,2$, we have $p\mid |K_i/H\cap K_i|$.
%Thus $\ZZ_p^2\leqslant K_1/(H\cap K_1)\times K_2/ (H\cap K_1) \lesssim G/H$ which is not possible. So $G_i$  ($i=1,2$) are abelian groups.
%
%
%Now suppose $G_1\cong G_2$ are non-cyclic metacyclic groups.
%Then there is a prime $p$ such that $p\mid (|H_1|,|K_1|)$.
%So $\ZZ_p^4\leqslant G$ which is not possible as $G$ is a metacyclic group. So $G_i$ $(i=1,2)$ are cyclic groups.\qed
%

\begin{lemma}[Kazarin \cite{Kazarin}]\label{T-2'}
 Let $T$ be a non-abelian simple group which has a Hall
$2'$-subgroup. Then $T=\PSL(2,p)$, where $p=2^e-1$ is a prime. Further, $T=GH$, where $G=\ZZ_{p}{:}\ZZ_{\frac{p-1}2}$ and $H=\D_{p+1}=\D_{2^e}.$
\end{lemma}

%Let $G$ be a metacyclic group, then $G$ has a beautiful "p-tower" structure \cite{meta-odd}.
%\begin{proposition}
%If $p_1<p_2<\cdots<p_k$ is the increasing sequence of all prime divisors of the order of $G$, then each Hall $\pi_i$-subgroup of $G$ is normal for $\pi_i=\{p_i,\cdots,p_k\}$. In particular, the Sylow subgroup for the largest prime divisor is normal.
%\end{proposition}

\section{Automorphisms of edge-transitive   metacirculants}

In this section, we study automorphism groups of  tetravalent edge-transitive   metacirculants. We first give a simple lemma on automorphism groups of tetravalent graphs of odd order.

\begin{lemma}\label{regular}
Let $\Ga$ be a tetravalent graph of odd order. Let $X\leqslant \Aut\Ga$ be vertex-transitive on $\Ga$. Suppose further that $X$ is of odd order. Then $X$ is regular.
\end{lemma}
\proof
Let $\a$ be a vertex.
Suppose that $X_\a\not=1$. Then $X_\a$ is a 3-group since $|G|$ is odd and $|\Ga(\a)|=4$, and 
$X_\a$ divides $\Ga(\a)$ into two orbits, one of which has length 3, and the  other has length 1. Since $\Ga$ is undirected graph, we have $\Ga$ has an orbital graph $\Sig$ which is of valency $1$. This is not possible, as  $\Sig$ is a graph of odd order. Hence $X_\a=1$, and $X$ is regular.\qed

%Given a group $G$ and a subset $S$ with $S=S^{-1}=\{s^{-1}\mid s\in S\}$, the Cayley graph 
%of $G$ with respect to $S$ is denoted by $\Cay(G,S)$.
%Then $\Ga(G,S)$ is connected if and only if $\l S\r=G$.
%A graph $\Ga=(V,E)$ is a Cayley graph of $G$ if and only if $\Aut\Ga$ contains a subgroup isomorphic to $G$ and transitive on $V$.
%\begin{lemma}\label{cayley1}
%Let  $\Ga=(V,E)$ be a connected tetravalent   metacirculant of odd order.
%Let $X=\Aut\Ga$, and $\a$ be vertex.  
%Then 
%\begin{itemize}
%\item [i)]$X_\a$ is a $2$-group. 
%\item[ii)] Let $G$ is a vertex-transitive group of odd order, then $G$ is regular and metacyclic.
%\item[iii)] $\Ga=\Cay(G,S)$ for some metacyclic group $G$.
%\end{itemize}
%\end{lemma}

\begin{lemma}\label{cayley}

Let $\Ga=(V,E)$ be a connected tetravalent   metacirculant of odd order. Then $\Ga=\Cay(G,S)$ for some metacyclic group $G$. Let $X\leqslant \Aut\Ga$ be such that $G\leqslant X$. Then each normal subgroup of $X$ of odd order is a subgroup of $G$.
\end{lemma}

\proof By the definition of   metacirculant, $\Aut\Ga$ has  transitive metacyclic subgroups. Let $G$ be such a subgroup of smallest order.
We claim that  $G$ is of odd order. Suppose $|G|$ is even. Since $G$ is metacyclic, we have $G$ is soluble. Then a Hall $2'$-subgroup $G_{2'}$ is of  odd order. Since $G$ is transitive, $X=GX_\a$. Moreover, $G_{2}\leqslant X_\a$. Thus $X=G_{2'}X_\a$, so $G_{2'}$ is transitive on $V.$
 Now $G_\pi < G$ which is contradict to our assumption. So $G$ is of odd order. By Lemma~\ref{regular}, $G$ is regular, so $\Ga=\Cay(G,S)$.

Let $M$ be a normal subgroup of $X$ of odd order. Let $Y:=MG$. Then $Y$ is transitive on $V\Ga$, and $|Y|$ is odd order. By Lemma~\ref{regular}, $Y$ is regular. So $M\leqslant  G$\qed

From now on, we suppose $\Ga$ is a connected edge-transitive tetravalent   metacirculant of odd order.

By the result of \cite{Li-s-arc}, there is no 4-arc transitive graph of valency at least 3 on odd number of vertices. The following proposition chartcterizes the vertex stabilisers of connected $s$-transitive tetravalent graphs (see \cite[Theorem 1.1]{Li-s-arc}, \cite[Lemma 2.5]{LLD}, and \cite[Proposition 2.16]{ZF1})

\begin{proposition}\label{stabilizer}
Let $\Ga$ be a connected $(X, s)$-transitive tetravalent graph of odd order, and let $X_\a$ be the stabilizer of a vertex $\a\in V$ in $X$. Then one of the following holds:
\begin{itemize}
\item[1)] $s=1$, $X_\a$ is a $2$-group;
\item[2)] $s=2,\A_4\leqslant X_\a\leqslant \S_4$;
\item[3)] $s=3$, $\A_4\times \ZZ_3\lhd X_\a\lhd \S_4\times\S_3$.
\end{itemize}
\end{proposition}
%

%\begin{lemma}\label{cayley}
%Let $G$ be a finite group of odd order, and let $\Ga=\Cay(G,S)$ be a connected $X$-vertex-transitive graph of valency $4.$ Assume that $M\lhd X\leqslant \Aut\Ga$ such that $M$ has odd order and $G\leqslant X$.  Then $M\leqslant G$.
%\end{lemma}
%\proof
%Let $Y:=MG$. Then $Y$ is transitive on $V\Ga$.
%Suppose $M\nleqslant G$. Then $|Y_\a|=|Y|/|G|=\frac{|M|}{|M\cap G|}.$
%Since $M$ is of odd order, and $\Ga$ is of valency $4$, we have  $Y_\a$ is a $3$-group.
%Thus $Y_\a$ divides $\Ga(\a)$ into two orbits, one of which has length 3, and the  other has length 1. It follows that $\Ga$ has an orbital graph of valency $3$. Let $\Sig$ be such an graph. Then $\Sig$ is a graph of odd order and out-valency $3$. Thus $\Sig$ is directed and hence $\val\Sig=6>4$ which is not possible. So $M\leqslant G.$\qed
%The {\it Frattini subgroup} $\Phi(G)$ of a group $G$  is the intersection of all maximal subgroups of $G$. Thus $\Phi(G)\char G$.
%
%Now we give our key lemma which shows that if $\Ga$ is $X$-edge-transitive then $\Ga$ is a Cayley graph of a split metacyclic group, and  either $X$ is insoluble, or $G$ is normal in $X$.
The following lemma characterises tetravalent  edge-transitive   metacirculants which has a quasiprimitive automorphism group.

\begin{lemma}\label{edge-transitive}
Let $\Ga$ be a connected tetravalent $X$-edge-transitive   metacirculant of odd order. Then
 $\Ga$ is a split metacirculant  of $G$. Suppose that $ X$ is quasiprimitive on $V$ and contains $G$. Then $X$ is almost simple.

% \begin{itemize}
							
%\item[(i)] $X$ is insoluble, and $(X,G)$ is listed in Table 2, or
%\item[(ii)] $G$ is normal in $X$ and $G\cong \ZZ_p$ or $\ZZ_p^2$ with $p$ a prime.
%\end{itemize}
\end{lemma}

\proof 
% Let $\Ga=(V,E)$ be a connected tetrevalent edge-transitive   metacirculant of odd order.
 By Lemma   ~\ref{cayley}, there is a metacyclic group $G$ such that $\Ga=\Cay(G,S)$.

 Since $X$ is quasiprimitive on $V$. 
 Let $M$ be a mininal normal subgroup of $X$. 
 Then $M=T^l$ where $T$ is a simple group.
  Suppose $X$ is affine.
 Then $M=\ZZ_p^d$  is transitive on $V$, and $M$ is of odd order. 
 By Lemma \ref{cayley}, $M\leqslant G$, so $d\leqslant 2$ as $G$ is metacyclic.
 Since $M$ is transitive on $V$, we have $G= M$, and $G$ is normal in $X$.
 Thus $X=G{:}\Aut(G,S)$. 
 Hence $\Aut(G,S)\leqslant \S_4$, as $\Aut(G,S)$ acts on $S$ faithfully.
 So $X_\a\leqslant \D_8$ as neither $\A_4$ nor $\S_4$ has a faithful representation of degree 2. 

Now suppose $T$ is insoluble.
We claim that $l=1$ and $X$ is almost simple. 
First of all, we notice that $X_\a$ is a $\{2,3\}$-group, so $|X|_{\{2,3\}'}=|G|_{\{2,3\}'}$.
Since $G$ is soluble, a Hall $\{2,3\}'$-subgroup $G_{\{2,3\}'}$ exists and $G_{\{2,3\}'}<G<X$, so $G_{\{2,3\}'}$ is a Hall $\{2,3\}'$-subgroup of $X$. 
Thus $X_{\{2,3\}'}$ is metacyclic.

Suppose $l> 2$. Let $p$ be a prime divisor of $|T|$ such that $p\neq 2,3$. Then a Sylow $p$-subgroup of $T^l$ is isomorphic to $T_p^l\leqslant X_p$, and  $X_p$ is a metacyclic group which is not possible.

Suppose $l=2$. 
Suppose $3\nmid |X_\a|$. Then $X_\a$ is a 2-group.
Thus $G_{2'}$ is a Hall $2'$-group of $X$. By Lemma~\ref{Hall}, $T$ has a Hall $2'$-subgroup $T_{\{2'\}}$.
Hence  $T\cong \PSL(2,p)$ with $T_{2'}=\ZZ_{p}{:}\ZZ_{\frac{p-1}2}$ by Lemma ~\ref{T-2'}. On the other hand, $T_{2'}^2\leqslant X_{2'}\cong G$. By Lemma~\ref{non-iso-direct}, $T_{2'}$ is a cyclic group. So $p=3$, and $T=\A_4$ is soluble which contradicts to $T$ is insoluble. 
Suppose $3\mid |X_\a|$.
By Lemma ~\ref{stabilizer}, $|X_\a|_2\leqslant 16$,
so $|T_2|\leqslant 4$ as $(T^2)_2\leqslant X_\a$. Check the order of non-abelian finite simple group, we have $T=\PSL(2,q)$ with $q\equiv \pm 3\pmod 8$ or $q=2^2$. Moreover, $T_{\{2,3\}'}^2\lesssim X_{\{2,3\}'}\cong G $. So $T_{\{2,3\}'}$ is a cyclic group. 

The order of $\PSL(2,q)$ is $q(q^2-1)/d=q(q+1)(q-1)/d$, where $d=(2,q-1).$
Suppose $p\geqslant 5$. Then the order of cyclic subgroup of $\PSL(2,q)$ is at most $(q+1)/2$ or $p$, if $q=p^f$ with $f\geqslant 2$, or  $q=p$, respectively, where $p$ is a prime. Suppose $3$ divides $q$. Then $2,3$ do not divide $(q+1)(q-1)/8$ as $q\equiv\pm 3\pmod 8$. Thus $(q+1)(q-1)/8=|T|_{\{2,3\}}'$. Hence there is a cyclic subgroup of order $(q+1)(q-1)/8$, so $(q+1)(q-1)/8\leqslant (q+1)/2$ or $q$. Since $3\mid q$ and $q=p^f>5$, we have the latter is not possible. So $(q+1)(q-1)/8\leqslant (q+1)/2$, thus $q\leqslant 5$ and $q=5$ contradicts to $3\mid q$. 
Suppose $3\mid p+1$. Similarly, $2,3$ do not divide $q(q-1)/8$. It follows that either $q(q-1)/8\leqslant q+1/2$ where $q=p^f$ with $f\geqslant 2$, or 
$q(q-1)/8\leqslant q$ where $q=p$. Thus $q\leqslant 7$, while $3\mid q+1$ so $q=5$. Suppose $3\mid q-1$. Using the same argument, we have $q=7$.
Notice that $\PSL(2,4)\cong \PSL(2,5)$, we have $T$ is $\PSL(2,5)\cong \A_5$, or $\PSL(2,7).$
For the former, $M=\A_5^2$, $X\leqslant M.\ZZ_2^2.\ZZ_2$ as $X$ is quasiprimitive. Hence the maximal metacyclic subgroup of odd order in $X$ is $\ZZ_5^2$. So $|V|=5^2$. On the other hand, $M_\a\lhd X_\a$ such that $X_\a$ satisfies Lemma~\ref{stabilizer}. Since the maximal subgroup of $M$ which is normal in $X_\a$ is $\A_4\times \ZZ_3$, we have $5^2=|V|=|M{:}M_\a|=4\times 5^2$ which is not possible. 
For the latter, $M=\PSL(2,7)^2$. Using the same argument we can show that this is not possible too.
Hence $l=1$ and $X$ is almost simple. \qed
%By \cite[Theorem 4.1]{LP}, $\A_{n-1}\leqslant X_\a\leqslant \S_{n-1}$, or subgroup $P_1$, $P_{n-1}$ of $\PSL(n,q)$, or $(X,G,X_\a)$ is given. Since for all these cases, $3$ divides $|X_\a|$. Thus $\Ga$ is $(X,s)$-arc-transitive with $s\geqslant 2$. By Lemma~\ref{stabilizer}, either $\A_4\leqslant X_\a\leqslant \S_4$ or $\A_4\times \ZZ_3\leqslant X_\a\leqslant \S_4\times \S_3$.
%Thus either $n=5$ and $X=\A_5,\S_5$, or $X$ are given in Tables 1 and 4. Checking the candidates given in Tables 1 and 4 in \cite{LP}, such that $G$ is of odd order and $X_\a$ satisfies the above conditions
%
%
%
%Thus $X$ is primitive  and arc-transitive.
%By \cite[Theorem 1.5]{Li-s-arc}, $(X,X_1,G)$  is listed in Table 1 and for each $G$ there is only one vertex-primitive arc-transitive graph of valency 4.\qed
%
%In the following, we suppose $X\leqslant \Aut\Ga$ is not quasiprimitive on $V.$ 
%

Quasiprimitive permutation groups which contain a transitive metacyclic subgroup
are determined in \cite{LP}. Then we have the following result.

\begin{lemma}\label{quasipri-gps}
Let $H$ be an almost simple quasiprimitive permutation group which contains
a transitive metacyclic subgroup $R$ of odd order. Then either $(H,H_\a)=(\A_n,\A_{n-1})$
or $(\S_n,\S_{n-1})$ with $n$ odd, or the triple $(H,R,H_{\a})$ lies in Table 2,
where $p$ is a prime, $G(q^n)$ is a transitive subgroup of
$\GammaL(1,q^n)=\ZZ_{q^n-1}{:}\ZZ_n$, and $P_1$ is a parabolic subgroup.
\end{lemma}
\begin{table}[htb]
%\caption{Table 2: quasiprimitive permutation groups of type \AS}
\[\begin{array}{lllll} \hline
\mbox{row} & H & R & H_{\a} & \mbox{conditions}\\ \hline

1 & \A_p.o & p{:}{p-1\over2} & \S_{p-2}\times o & o\le 2 \\

2 & \PSL(n,q).o & G(q^n).o_1 & P_1.o_2 & q=p^f,~o_1o_2\cong o\leq f.(n,q-1)\\

3   &\PSL(2,p).o & p{:}{p-1\over2} & \D_{(p+1)o}& o\le 2,\ p\equiv 3\ (\mod 4)\\

4 & \PSL(2,11).o & 11{:}5 & \A_4.o & o\leq 2 \\

5  & \PSL(2,11) & 11 & \A_5 & \\

6 & \PSL(2,29) & 29{:}7 & \A_5 & \\

7 & \PSL(2,p) & p{:}{p-1\over2} & \A_5 & p=11,19,59 \\

8 & \PSL(2,23) & 23{:}11 & \S_4 & \\

9  & \PSL(5,2) & 31{:}5 & 2^6:(\S_3\times \PSL(3,2))& \\

10 &\PSU(3,8).3^2.o & 3\times19{:}9 & (2^{3+6}{:}63{:}3).o &
o\le2\\

11 & \PSU(4,2).o & 9{:}3 & 2^4{:}\A_5.o & o\le 4\\

12 &\M_{11} & 11,\ 11{:}5 & \M_{10},\ \M_9.2 & \\

13  &\M_{23} & 23,\ 23{:}11&\M_{22},\ \M_{21}.2,\ 2^4.\A_7 & \\ \hline

\end{array}\]
\centerline{Table 2: quasiprimitive permutation groups of type \AS}

\nobreak
%\centerline{\bf Table~A} \vskip0.08in
\end{table}

\begin{lemma}\label{AS}
Suppose $\Ga=\Cay(G,S)$ is a connected $X$-edge-transitive tetravalent graph of odd order such that $G\leqslant X$ is quasiprimitive group of almost simple. Then $(X,G,X_\a)$ lies in Table 3.
\end{lemma}

\proof First of all, since $X_\a$ satisfies Lemma~\ref{stabilizer} and $G$ is of odd order, by checking the candidates in Lemma\ref{quasipri-gps}, we have $(X,X_\a)$ can only appear as rows 3,4,8 or $(\A_5,\A_4)$, $(\S_5,\S_4)$.
Suppose $(X,X_\a)=(\PSL(2,p).o,\D_{p+1}.o)$. Then by Lemma~\ref{stabilizer}, $\D_{p+1}.o$ is a 2-group. So $X$ is primitive except $p=7$. Since $X$ are primitive for all the other five cases, we have either $X$ is primitive, or $X=\PSL(2,7)$ and $X_\a=\D_8$. Suppose $X$ is primitive. Then by  \cite[Theorem 1.5]{Li-s-arc}, $(X,X_1,G)$  is listed in rows 1,3 and 4 of Table 1 and for each $G$ there is only one vertex-primitive arc-transitive graph of valency 4. Suppose $X=\PSL(2,7)$ with $X_\a=\D_8$. Then $X_{\a\b}=\ZZ_2$, and $\Aut\Ga=\PGL(2,7)$ which is primitive too. So this gives the candidate of  row 2 of Table 1.\qed

\begin{lemma}\label{normal}
Let $\Ga=\Cay(G,S)$ be a connected $X$-edge-transitive tetravalent graph of odd order. Suppose $G\leqslant X$ has a non-trivial intransitive normal subgroup $M$.
Then $G$ is normal in $X$, and $X_\a\leqslant \D_8$.
\end{lemma}

\proof
Let $M$ be a minimal normal subgroup of $X$ which is intransitive. 
Let $\BB$ be the set of $M$-orbits. Let $K$ be the kernel of $X$ acting on $\BB$.

{\bf Case 1. } Suppose $K_\a\neq 1$. 
We claim that $G\lhd X$. 
Suppose $G$ is not normal in $ X$. 
By \cite[Lemma 5.6]{LLZ}, $G=\ZZ_p^2{:}\ZZ_m$ with $m\geqslant3$ 
and $X=\ZZ_p^2{:}(\ZZ_2^l{:}\ZZ_m)$ or $\ZZ_p^2{:}(\ZZ_2^l{:}\ZZ_m).\ZZ_2$ such that
$\ZZ_p^2$ is the only minimal normal subgroup of $X$, where $l\geqslant2$.
It follows that $\ZZ_2^l{:}\ZZ_m$ has a faithful representation on $\ZZ_p^2$, which is not possible. 
Thus $l=1$ and
%Since $K_\a\neq 1$ and $\Ga$ is $X$-edge-transitive, we have $\Ga_M$ is a cycle. 
%Suppose $\a\in B_0$. 
%Then there are exactly two vertices, say $\b,\g\in B_1$ who are adjacent to $\a$ in $\Ga$.
%Thus the  induced graph $[B_0\cup B_1 ]$ is a union of cycles, say $\cup[B_{0s},B_{1s}]$, where $[B_{0s},B_{1s}]$ are  cycles of length $2|B_0|/s$.
%So $2|B_0|/s>4$.  
%Since $\Ga$ is edge-transitive and connected, for every cycle $[B_{ij},B_{i+1j}]$, where $B_{ij}\subset B_i$ and $B_{i+1j}\subset B_{i+1,}$there is a "cycle path" $(B_{01},B_{11},\cdots B_{ij},B_{i+1j})$. 
%Let $ g\in \K_\a$ such that $g$ fixes $\b$ and $\g$. 
%Then $g$ fixes three vertices $\a,\b$ and $\g$ in $[B_{01}\cup B_{11}]$. 
%So $g$ fixes all vertices in $B_{01}\cup B_{11}$.
%It follows that $g$ fixes all vertices in $B_{11}\cup B_{21} $ and so on.
%Thus $g$ fixes $V$. Hence $g=1$.
%So $K_\a$ acts on $\{\g,\b\}$ faithfully and 
$K_\a\cong \ZZ_2.$
So $K=M:\ZZ_2$. 
By Lemma   ~\ref{cayley}, $M<G$. Let $\ov G:=GK/K\cong G/ K\cap G $, let $\ov X=X/K$. We have $\ov G\leqslant \ov X$ is transitive on $\Ga_M$. 
Since $G$ is of odd order, $\Ga_M$ is a cycle, we have $\ov G$ is regular on $V\Ga_M$, so $\ov G\lhd\ov X$. 
Then $K.\ov G\lhd X$. 
Let $L=K.\ov G=K.(G/G\cap K)$.
Since $K=(K\cap G){:}\ZZ_2$, we have  $L=2|V|$ is soluble. 
So $G=L_{2'}$ has index 2 in $L$.
Thus $G\char L\lhd X.$
Hence $G\lhd X$.

{\bf Case 2.} Suppose $K_\a=1.$ We will show that $G\lhd X$.
Suppose for each normal subgroup $N$ of $X$, $\Ga$ is a normal cover of $\Ga_N$. Then let $N$ be a maximal normal subgroup of $X$ which is intransitive. It follows that $X/N$ is quasiprimitive on $V\Ga_N$. 
Thus by Lemma   ~\ref{cayley}, we have $N<G$.
Suppose further $N$ is a minimal normal subgroup of $X$.  
Then $N=\ZZ_p^d$ with $d\leqslant 2$ as $N<G$. 
Moreover, $X/N$ is quasiprimitive on $\Ga_N$. 
So either $X/N$ is an affine group such that $G/N\lhd X/N$ or $X/N$ is listed in Table 1. For the former, since $G/N\lhd X/N$, we have $G\lhd X$.
For the latter, $N=\ZZ_p^d$ with $d\leqslant 2$ and $p$ an odd prime. By  \cite{Low-dimension}, the groups listed in Table 1 are not contained in $\Aut(N)=\GL(2,p)$ or $\GL(1,p)$. Further, as $p$ is an odd prime, and the Schur Multipliers of $\A_5,\PSL(2,7),\PSL(2,11)$ and $\PSL(2,23)$ are 2, we have $X=N\times (X/N)$. Hence $X$ has a normal subgroup $H$ which is isomorphic to the socle of $X/N$. Further,  $H$ is intransitive on $V$. So $H$ has at least three orbits on $V$, but $3\mid |H_\a|$ which is not possible. 
Suppose $N$ is not minimal normal subgroup of $X$.
Let $N'\lhd N$ be a  second maximal normal subgroup of $X$ which is intransitive. 
Then  $N/N'$ is a minimal normal subgroup of $  X$. Then $\Ga$ is a normal cover of $\Ga_{N'}.$
Let $\ov N:=N/N'$, and $\ov X:=X/N'$. 
Then $\Ga_{N'}$ is $\ov X$-edge-transitive with $\ov G:=G/N'\leqslant X/N'$ is regular on $V\Ga_{N'}$, and $\Ga_{N'}$ is a normal cover of $\Ga_N$. 
Moreover, $\ov X/\ov N\cong X/N$ is quasiprimitive on $V\Ga_N$. 
By the same argument as above, we have this is not possible.

Now suppose there is a normal subgroup $L$ of $X$, such that $\Ga$ is not a normal cover of 
$\Ga_L$.
Since $\Ga$ is a normal cover of $\Ga_M$, we may suppose there are normal subgroups $N,L$ of $X$ such that $\Ga$ is a normal cover of $\Ga_N$, $L>N$ such that $L/N$  is a minimal normal subgroup of $X/N$ and $\Ga$ is not a normal cover of $\Ga_L$.
By Lemma \ref{cayley}, $N<G$.
Now $\ov X:=X/N$ is edge-transitive on $\Ga_N$ and $\ov G:=G/N\leq X/N$ is regular on $V\Ga_N$. Further $\ov L:=L/N$ is a minimal normal subgroup of $\ov X=X/N$ and $\Ga_N$ is not a normal cover of $\Ga_L$. By Case 1. we have $G/N\lhd X/N$, so $G\lhd X$.

Thus $X=G{:}\Aut(G,S)$. 
 Since $\Aut(G,S)$ acts on $S$ faithfully, we have $\Aut(G,S)\leqslant \S_4.$ 
 So $X_\a\leqslant \D_8$ as neither $\A_4$ nor $\S_4$ has a representation of degree 2. \qed

Now we study normal edge-transitive metacirculants. 

\begin{lemma}\label{normal-meta}
Let $\Ga=\Cay(G,S)$ be a connected $X$-edge-transitive metacirculants of odd order such that $G\lhd X$. Then one of the following holds:
\begin{itemize}
\item[1)] $\Ga$ is  an $X$-half-transitive metacirculant of non-abelian metacyclic group, and  $X=G{:}\ZZ_2$;

\item[2)] $\Ga$ is an $X$-arc-transitive  metacirculant of an abelian metacyclic group, and $X=G{:}X_\a$, where $\ZZ_2<X_\a\leqslant\D_8$.
\end{itemize}
\end{lemma}

\proof
Let $\Ga=\Cay(G,S)$, where $S=\{x,x^{-1},y,y^{-1}\}$ is a generating subset of $G$.
By Lemma ~\ref{normal}, $X\leqslant G{:}\Aut(G,S)$ with $\Aut(G,S)$.
Let $\s\in X_\a$ where $\a\in V$.
Suppose $x^\s=x^{-1},$ $y^{\s}=y^{-1}$. 
Then $o(\s)=2$. 
Suppose further $\s$ fixes $z=x^iy^j$ for some integers $i,j$.
Then \[x^iy^j=(x^iy^j)^\s=x^{-i}y^{-j},\]
and so $x^{2i}=y^{-2j}$. Since $o(x)=o(y)$ are odd, we have $x^i=y^{-j}$.
Then $z=x^iy^j=1$, that is, $\s$ fixes no non-identity of $G$. 
By \cite[Exercise 1.50]{R}, $G$ is abelian. Since all edge-transitive Cayley graphs of abelian groups are arc-transitive, we have $\Aut(G,S)>\ZZ_2$.
Further, by Lemma  ~\ref{edge-transitive}, $\Aut(G,S)=\D_8.$

Suppose $G$ is non-abelian.
Suppose there is an element $\t\in\Aut(G,S)$ such that $x^\t=x^{-1}$, or $y^\t=y^{-1}$. 
Without loss of generality, we suppose $x^\t=x^{-1}$.
Since $G$ is non-abelian, $y^\t\neq y^{-1}$. 
So  $y^\t=y$.
As $\Ga$ is edge-transitive, $\Aut(G,S)$ has at most two orbits  on $S$, 
So there is another element say $\g\in \Aut(G,S)$ such that $y^\g=y^{-1}$, and $x^\t=x$.
Thus $x^{\g\t}=x^{-1}$, $y^{\g\t}=y^{-1}$, and $o(\g\t)=2$ which contradicts to $G$ is non-abelian.
Without loss of generality, we suppose $x^\t=y$ and suppose further that  $y^\t\neq x$. Then $y^\t=x^{-1}$.
Hence $x^{\t^2}=x^{-1}$ which  is not possible too.
So $y^\t=x$. Thus $\Aut(G,S)\cong \ZZ_2.$\qed

%\begin{example}
%{\rm
%\begin{itemize}
%\item[1.] Let $G_1=\l a,b\mid a^{5}=b^5=1,ab=ba\r$. Let $S_1=\{a,b,a^{-1},b^{-1}\}$.
%Let $\Ga_1=\Cay(G_1,S_1)$. 
%Then $\Ga_1$ is arc-transitive and $\Aut\Ga_1=\ZZ_5^2{:}\D_8$.
%
%\item[2.] Let $G_2=\l a,b\mid a^{5^2}=b^5=1,ab=ba\r$, $i=1,2$ Let $S_2=\{ab,a,a^{-1}b^{-1},a^{-1}\}$.
%Let $\Ga_2=\Cay(G_2,S_2)$. 
%Then $\Ga_2$ is arc-transitive and $\Aut\Ga_2=(\ZZ_{5^2}{:}\ZZ_5){:}\ZZ_2^2.$
%
%
%\item[3.] Let $G_3=\l a\r\cong \ZZ_p$ such that $p$ is a prime and $p\equiv 1\pmod 4$ . Let $S_3=\{a,a^{-1},a^2,a^{-2}\}$. Let $\Ga_3=\Cay(G_3,S_3)$Then $\Aut\Ga=\ZZ_p{:}\ZZ_4$.
%
%\item[4.] Let $G_4=\l a,b\mid a^p=b^q\r\leqslant \AGL(1,p)$ where $p$ is an odd prime and $q$ is odd. Let $S_4=\{b,ab,b^{-1},b^{-1}a^{-1}\}$. Let $\Ga_4=\Cay(G_4,S_4) $. Then $\Aut\Ga=\ZZ_p{:}\ZZ_{2q}.$ 
%\end{itemize}}
%\end{example}
%

%\vskip0.1in 
%
%{\bf Proof of Theorem  ~\ref{thm-1}:} By Lemmas  ~\ref{edge-transitive} and   ~\ref{normal-meta}, Theorem  ~\ref{thm-1} holds.

\section{Edge-transitive tetravalent Sylow-circulants}

As an application of Lemma ~\ref{normal-meta}, in this section we study  Sylow-circulants of odd order and determine the number of non-isomorphic edge-transitive tetravalent  Sylow-circulants.

 Recall that a group is called Sylow-cyclic, if all Slow subgroups of it are cyclic.
 Let $G=H\times K=\l c\r\times(\l a,b\r)=\l a\r{:}\l bc\r$ be a Sylow-cyclic group, where $K=\l a,b\mid a^m=b^n=1,a^b=a^r\r$, $r\not\equiv 1\pmod m, \ r^n\equiv 1\pmod m$, such that no Sylow-$p$ subgroup of $\l b\r$ lies in the centre of $G$. We first study the edge-transitive tetravalent Cayley graph of $G$ such that $H=1$.
%\[ a^ub^v=b^va^{r^vu}\]
%\[(a^{u_1}b^{v_1})(a^{u_2}b^{v_2})=b^{v_1+v_2}a^{u_1r^{v_1+v_2}+u_2r^{v_2}}\]
%\[(a^ub^v)^k=b^{kv}a^{u(r^{kv}+r^{(k-1)v+\cdots+r^v})}=b^{kv}a^{u[r^v]_k}\]
%\vskip0.1in
%
%{\bf Remark:} By equation \ref{eqn:3}, we have $o(a^ib^j)=o(a^{i[r^j]_l})\cdot l$, where $o(b^j)=l$.
%
%\vskip0.1in
%
%As proved in Lemma \ref{edge-transitive}, if $\Ga=\Cay(G,S)$ is $X$-edge-transitive  such that $G\leqslant X$, where $X$ is soluble, then $X=G{:}\Aut(G,S)$, where $\Aut(G,S)\leqslant\Aut(G)$.
% Thus the group $\Aut(G)$ is very important. 
By Lemma ~\ref{normal-meta}, we know all edge-transitive tetravalent graphs $\Cay(G,S)$  are normal Cayley graphs.
So $\Aut\Ga=G{:}\Aut(G,S)$.
Thus to get a better understand of the automorphism groups and these graphs, we will study the automorphism group $\Aut(G)$ of $G$, then determine $\Aut(G,S)$. 
\vskip0.1in

{\bf Remark:} Let $G=\ZZ_m{:}\ZZ_n$ be a Sylow-cyclic metacyclic group. Then $(m,n)=1$.
\vskip0.1in

Now $G=\l a,b\mid a^m=b^n=1, a^b=a^r\r$ is a non-abelian metacyclic group. Then  $r\not\equiv  1\pmod m$ and $r^n\equiv 1\pmod m$. 
%For convenience, denote $[r]_s=r^{s}+r^{s-1}+\cdots+r$, where $s$ is an positive integer.
%It is easily seen that $[r]_s$ divides $[r]_{st}$ for positive integers $s$ and $t$.
%
By calculating, the following equations hold.
\begin{equation}\label{eqn:1}
a^ub^v=b^va^{r^vu}
\end{equation}
\begin{equation}\label{eqn:2}
(a^{u_1}b^{v_1})(a^{u_2}b^{v_2})=b^{v_1+v_2}a^{u_1r^{v_1+v_2}+u_2r^{v_2}}
\end{equation}\begin{equation}\label{eqn:3}
(a^ub^v)^k=b^{kv}a^{u(r^{kv}+r^{(k-1)v+\cdots+r^v})}=b^{kv}a^{ur^v\frac{r^{vk}-1}{r^v-1}}
\end{equation}

\begin{lemma}\label{order}
Let $G=\l a, b\mid a^m=b^n=1,a^b=a^r\r$ be a non-abelian Sylow-cyclic group. Suppose further that $G$ is not a product of its any two non-trivial subgroups. Then  $o(a^ib^j)=o(b^j)=n$ for $0\leqslant i\leqslant m$ and $(j,n)=1$.
\end{lemma}

\proof We first prove that $(r-1,m)=1$ and $\l a\r\cap \Z(G)=1$.
%Suppose $a^l\in \Z(G)$ such that $p\mid o(a^l)$, where $p$ is a prime. Then $(a^l)^b=a^{lr}=a^l$. So $l(r-1)\equiv 0\pmod m$. Thus $p\mid (r-1,m).$ 
Let 
$p^\a$ be the largest $p$-power divisor of $(r-1,m)$. Then  $r\equiv p^\a x+1\pmod m$ where $(x,p)=1.$ 
Since $r^n\equiv 1\pmod m$, we have $(r^n-1)\equiv (p^\a x)^n+n(p^\a x)^{n-1}+\cdots+n(p^\a x)\pmod m$. 
Thus $p^\a$ is the largest prime divisor of $m$ as $(n,m)=1$.  
It follows that if $\a\neq 0$, then  $(a^{m/p^\a})^b=a^{rm/p^\a }=a^{m/p^\a}$, that is $\l a\r_p\subset \Z(G)$. So $G=\l a\r_p\times G'$ which contradicts to our assumption. So $\a=1$ and $(r-1,m)=1$. Suppose $a^l\in\Z(G)$.
Then $(a^l)^b=a^{rl}=a^l$. Hence $l\equiv 0\pmod m$ as $(r-1,m)=1$. So $\l a\r\cap\Z(G)=1$.
Moreover, $(a^ib)^n=b^na^{ir\frac{r^n-1}{r-1}}=a^{ir\frac{r^n-1}{r-1}}.$ Since $(r-1,m)=1$, we have $a^{ir\frac{r^n-1}{r-1}}=1$, that is $(a^ib)^n=1$. So $o(a^ib)=n$ for any $i.$

Now suppose $p\mid (r^j-1,m)$. Then $(a^{m/p})^{b^j}=a^{r^j m/p}=a^{m/p}.$ Since $(j,n)=1$, we have $a^{m/p}\in\Z(G)$ which is not possible. So $(r^j-1,m)=1$ and then $o(a^ib^j)=n.$\qed

The following lemma give some properties of $\Aut(G)$.

 \begin{lemma}\label{auto-G}

Let $G=\l a,b\mid a^m=b^n=1, a^b=a^r\r$ such that $r\neq 1$ and $r^n\equiv 1\pmod m$.  Let $n_0$ be the smallest integer such that $r^{n_0}\equiv 1\pmod m$. 
Suppose further that $G$ is not a product of its any two non-trivial subgroups. 
Then $\s\in \Aut(G)$ if and only if:

\[\begin{array}{lll}
\s:&a \to a^s& \mbox{ where } 1\leqslant s<m, \  (m,s)=1\\
& b \to a^tb^{1+ln_0} & \mbox{ where } 0\leq l<n/n_0, 1\leqslant t\leqslant m.
\end{array}\]

In particular, $b^j$ is conjugate to $a^tb^{j+ln_0}$ under $\Aut(G)$ for $(j,n=1)$, and $b^i$ not conjugate to $b^{-i}$ under $\Aut(G)$ for any $b^i\not \in \Z(G)$.
\end{lemma}

\proof 
Let $\s$ be an automorphism of $G$. 
By assumption, $\l a \r\char G$, so $a^\s=a^s$ for some integer $s$ such that $o(a^s)=m$. Hence $(s,m)=1$. 
Now suppose $b^\s=a^tb^u$. Then 
\[
\begin{array}{lll}
&(ab)^\s&=a^\s b^\s=a^sa^tb^u=a^ta^sb^u\\
=&(ba^r)^\s&=a^tb^ua^{rs}.
\end{array}\]

Hence $a^sb^u=b^ua^{rs},$ that is, $b^{-u}a^sb^u=a^{rs}$. So $b^{-(u-1)}a^sb^{u-1}=b\big(b^{-u}a^sb^u\big)b^{-1}=ba^{rs}b^{-1}=a^s.$ 
Thus $b^{u-1}\in \Z(G)$, as $(s,m)=1$. Since $n_0$ is the smallest integer such that $r^{n_0}\equiv 1\pmod m$, we have $\Z(G)\cap\l b\r=\l b^{n_0}\r$. 
Thus $u=1+ln_0$ for some integer $0\leq l<n/n_0$.

 On the other hand, Let $\s$ be defined as in the lemma. Since $G=\l a^s, a^tb^{1+ln_0}\mid o(a^s)=m,\, o(a^tb^{1+ln_0})=n,(a^s)^{a^tb^{1+ln_0}}=(a^s)^r\r$, we have $\s$ is an automorphism of $G$. Similarly, suppose $a^\t=a^s,$ where $(s,m)=1$, and $(b^j)^\t=a^tb^{j+ln_0}$. Then $G=\l a^s, a^tb^{j+ln_0}\mid o(a^s)=m,\, o(a^tb^{1+ln_0})=n,(a^s)^{a^tb^{j+ln_0}}=(a^s)^{r^j}\r$. Then $\t$ is an automorphism of $G$, so $b^j$ is conjugate to $a^tb^{j+ln_0}$ under $\Aut(G).$
 
 Suppose $(b^i)^\s=b^{-i}$. 
Then $b^{i+lin_0}=b^{-i}$. 
Thus $2i+lin_0\equiv 0\pmod n$.
 Since $n_0\mid n$, we have $n_0\mid i$, and $b^i\in \Z(G)$. 
 Thus $b^i$ is not conjugate to $b^{-i}$ under $\Aut(G)$ for any $b^i\not \in\Z(G)$. So the lemma holds. \qed
  
 Let $\Ga=\Cay(G,S)$ be a connected edge-transitive metacirculants with $S=\{x,y,x^{-1},y^{-1}\}$. By Lemma  ~\ref{normal-meta}, we have $\Aut\Ga=G{:}\Aut(G,S)$ where $\Aut(G,S)=\ZZ_2$. Thus the involutions in $\Aut(G)$ play an important role in study edge-transitive metacirculants. Then following lemma determines the involution automorphisms of $G$.
 
Let $\Z$ be the centre of $G$. Let $n_0$ be the smallest integer such that $r^{n_0}\equiv 1\pmod m$. Then $a^{b^{n_0}}=a^{r^{n_0}}=a$. So $b^{n_0}\in \Z$.

\begin{lemma}\label{involution-auto}
Let $\s$ be an involution automorphism of $G$. Then 
\[\begin{array}{lll}
\s&:& b\to a^ib \qquad\mbox{  where }o(a^ib)=n\\
&&a \to a^s \qquad\mbox{  where } s^2\equiv 1\pmod m.
\end{array}\]

\end{lemma}
\proof Suppose $b^\s=a^ib^{1+ln_0}$, $a^\s=a^s$. 
Then 
\[b^{\s^2}=(a^ib^{1+ln_0})^\s=a^{is}(a^ib^{1+ln_0})^{(1+ln_0)}=b^{(1+ln_0)^2}a^{x}, \mbox{ as }[a,b]\in \l a\r.\]
Since $o(\s)=2$, $b^{(1+ln_0)^2}=b,$ that is, $(1+ln_0)^2\equiv 1\pmod n$. 
So $ln_0(ln_0+2)\equiv 0\pmod n$.
We know that $\l b^{n_0}\r\leqslant \Z$. Let $p$ be a prime divisor of $o(b^{n_0})$ and suppose $(p,n_0)=1$. 
Let $P$ be a Sylow-$p$ subgroup of $\l b^{n_0}\r$.
Since, $(p,n_0)=1$, we have $P$ is a Sylow-$p$ subgroup of $\l b\r$. Hence there is a Hall-$p'$ subgroup $B_{p'}$ of $\l b\r$ such that $\l b\r=P{:}B_{p'}$. 
While $P\leqslant \Z$, we $\l b\r=P\times B_{p'}$.
It follows that $G=P\times(\l a\r{:}B_{p'})$ which contradicts to our assumption.
 Thus $p$ divides $n_0$, so $(p,ln_0+2)=1$.
Hence $(ln_0+2,n)=1$, and $ln_0\equiv 0\pmod n.$\qed

\begin{lemma}\label{abj in S}
$S$ has the form $\{a^{i_1}b^j,a^{i_2}b^j,(a^{i_1}b^j)^{-1},(a^{i_2}b^j)^{-1}\}$ for some $1\leqslant i_1<i_2\leqslant m$, where $(j,n)=1$.
\end{lemma}

\proof  Since $\l S\r=G$, $\l \ov S\r=\ov G=G/\l a\r= \l \ov b\r\cong\ZZ_n$ where $\ov S,\ \ov b$ are the image of $S,\ b$ under the map from $G$ to $G/\l a\r,$ respectively. 
Thus there is an element, say $x$, with the form $a^{i_1}b^j$ where $(j,n)=1$. 

Further, by Lemma  ~\ref{normal-meta}, $\Aut(G,S)=\ZZ_2$. Let $\s$ be an involution in $\Aut(G,S)$. By Lemma  ~\ref{auto-G}, $\s$ can not maps $x$ to $x^{-1}$. Thus $\s$ map $x$ to $y$ (or $y^{-1}$). By Lemma  ~\ref{involution-auto},  $y=a^{i_2}b^j$.\qed

%\begin{lemma}\label{order}  $r^j\equiv r\pmod m$ if $(j,n)=1$.
%\end{lemma}
%
%\proof It is equivalent to prove $r^j-1\equiv r-1\pmod m$.
%Since $r^n\equiv 1\pmod m,$  we may suppose $r-1\equiv p^xy \pmod m$ such that $(y,m)=1$. Thus $r^j-1\equiv (p^xy+1)^j-1\pmod m$. Hence $p^x\mid r^j-1\pmod m$. It follows that $(r-1\pmod m) \mid (r^j-1\pmod m).$
%
%On the other hand, as $(j,n)=1$, there is an integer $i$ such that $(r^j)^i=r\pmod m$ where $(i,n)=1$. So $(r^j-1\pmod m) \mid (r-1\pmod m)$. Thus $r^j-1\equiv r-1\pmod m.$\qed
%

The following lemma gives a sufficient and necessary condition of $\l S\r=G$.
\begin{lemma}\label{generators}
Let $G=\l a\r{:}\l b\r$ as defined above, and let $S=\{a^{i_1}b^j,a^{i_2}b^{j}\}$,  $1\leqslant i_1,i_2\leqslant m$, $1\leqslant j\leqslant n$.
Then $\l S\r=G$ if and only if $(j,n)=1$,  $\Big( i_2-i_1,i_1[r]_n \pmod m\Big)=1$.
\end{lemma}

\proof Suppose $(j,n)=1$ and $(i_2-i_1,i_1[r]_n)=1$. 
Then $a^{i_2-i_1}=a^{i_2}b^j\cdot (a^{i_1}b^j)^{-1}\in \l S\r$, and $(a^{i_1}b^j)^n=a^{i_1[r^j]_n}=a^{i_1[r]_n}\in S$.
Since $(i_2-i_1,i_1[r]_n)=1$, we have $\l a\r\subset \l S\r$. 
Then $b^j\in \l S\r$.
Hence $\l b\r\subset \l S\r$, as $(j,n)=1$.

Now suppose $\l S\r=G$. 
Let $\phi$ be the natural map from $G$ to $G/\l a\r$. 
Let $\ov S$, $\ov b$ be the image of $S$ under $\phi$. 
Then $\l \ov S\r=\l \ov b^j\r$.
Since $\l S\r=G$, $\l \ov S\r=\l \ov b^j\r=\ov G\cong\ZZ_n$. 
Hence $(j,n)=1$.

Now $G=\l S\r=\l a^{i_1}b^j,a^{i_2-i_1}\r$. For convenience, denote $i=i_2-i_1$. 
 Then $\l a^i\r\lhd \l S\r.$ 
Thus every elements in $G$ can be written as $a^{is}(a^{i_1}b^j)^t$, where $s,t$ are integers. 
We may suppose $a=a^{is'}(a^{i_1}b^j)^{nt'}$. 
For convenience, let $x=(a^{i_1}b^j)^n=a^{i_1[r^j]_n}=a^{l}$.
Then $a=a^{is'+lt'}$
So $\Big(i,i_1[r^j]_n\pmod m\Big)=\Big(i,i_1[r]_n\pmod m\Big)=1$, that is, $\Big(i_2-i_1,i_1[r]_n\pmod m\Big)=1$.\qed

By Lemma  ~\ref{abj in S}, every edge-transitive tetravalent Cayley graph of non-abelian metacyclic group is isomorphic to $\Cay(G,S)$, where $S=\{a^{i_1}b^j,a^{i_2}b^j,(a^{i_1}b^j)^{-1},(a^{i_2}b^j)^{-1}\}.$ 
Moreover, as $\Aut\Ga=G{:}\ZZ_2$, $(|G|,|\Aut\Ga_\a|)=1$.
By Proposition  ~\ref{CI-criterion}, $\Ga$ is a CI-graph.

In the following, we will study the isomorphism between edge-transitive metacirculants.
Thus it is sufficient to determine the conjugation of different generationg set $S$.

\begin{lemma}\label{edge-trans-classification}
 Let $\Ga$ be a connected normal edge-transitive tetravalent graph. Then $\Ga\cong\Cay(G,S_j)$ with  $S_j=\{b^j,ab^{j},b^{-j}, (ab^{j})^{-1}\}$, where $(j,n)=1$ and $1\leqslant j<n_0$.

\end{lemma}

\proof 
Let $\Ga=\Cay(G,S)$. Since $\l S\r=G$, $\l \ov S\r=\ov G=G/\l a\r= \l \ov b\r\cong\ZZ_n$ where $\ov S,\ \ov b$ are the image of $S,\ b$ under the map from $G$ to $G/\l a\r,$ respectively. 
Thus there is an element, say $x$, with the form $a^{i_1}b^j$ where $(j,n)=1$. By Lemma ~\ref{auto-G}, $a^tb^{j+ln_0} $ is conjugated to $b^j$ under $\Aut(G).$ Thus $S$ is conjugate to $S_j'$, where $b^j\in S_j'$. Hence $\Ga\cong \Cay(G,S_j')$.
Further, by Lemma  ~\ref{normal-meta}, $\Aut(G,S_j')=\ZZ_2$. Let $\s$ be an involution in $\Aut(G,S_j')$. By Lemma  ~\ref{involution-auto},  suppose $\s$ map $b^j$ to $a^ib^j$. Then $S_j'=\{b^j,a^ib^j,b^{-j},(a^ib^j)^{-1}\}.$  It follows that $(i,m)=1$, as $\l S_j'\r=G$. Since $(i,m)=1$, we have there is an integer $s$ such that $is\equiv 1\pmod m.$
It follows that $(s,m)=1$.
Moreover, $(t_2,m)=1$, so $(st_2,m)=1$.
Let \[\t:\ b^j\to ab^j,\ a\to a^{-s}. \]
Then $S_j'^\t=S_j=\{b^j,ab^j,(b^j)^{-1},(ab^j)^{-1}\}$. 
Hence $\Ga_1\cong\Ga_2$. Thus $\Ga\cong\Cay(G,S_j)$ where $S_j=\{b^j,ab^j,b^{-j},(ab^j)^{-1}\}.$ The lemma holds.\qed

{\bf Proof of Theorem  ~\ref{thm-2}:}
Let $G=H\times K=\l c\r\times(\l a,b\r)=\l a\r{:}\l bc\r$, where $K=\l a,b\mid a^m=b^n=1,a^b=a^r\r$, $r^n\equiv 1\pmod m$, such that no Sylow-$p$ subgroup of $\l b\r$ lies in the centre of $G$. Suppose first that $H=1$. Then by Lemma  ~\ref{edge-transitive}, the theorem holds.

Now we suppose  $H\neq 1$.
Since $G$ is a metacyclic group, we have $(o(c),n)=1$.
Thus $H\char G$ and $K\char G$.
Let $S=\{x,y,x^{-1},y^{-1}\}$, where $x=x_1x_2$ and $y=y_1y_2$ with $x_1,y_1\in H$, $x_2,y_2\in K$.

As $\Ga$ is edge-transitive, by Theorem   ~\ref{thm-1}, there is an involution automorphism $\s$ such that $x^\s=y$. Thus $x_i^\s=y_i$, $i=1,2$.
Hence  $x_1\neq 1$ as $\l S\r=G$, and $\Cay(K,S_2)$ is a tetravalent edge-transitive metacirculant where $S_2=\{x_2,y_2,x_2^{-1},y_1^{-1}\}$. Thus by Lemma ~\ref{edge-transitive}, $S_2$ satisfies Lemma \ref{edge-transitive}.
Now consider $x_1,y_1$, as $\l S\r=G$, we have $\l x_1,y_1\r=H$.
Suppose $x_1=c^s$. Then $y_1=c^{st}$ with $t^2\equiv 1\pmod l.$
Thus  $(s,l)=1$  as $\l x_1,y_1\r=H$.
We may suppose $\{x_1,y_1\}=\{c,c^t\}$ as $c$ is conjugate to $c_s$ where $(s,l)=1$.
Then $\{x_1,y_1\}^\s=\{c,c^{-1}\}$, where $c^\s=c^t$. Thus by Proposition   ~\ref{CI-criterion}, the Theorem holds.

\end{document}